\documentclass[hyper,11pt]{JHEP3}

\usepackage{xy,amssymb,amsmath,diagrams,times,theorem,longtable}

\xyoption{poly}
\xyoption{arc}
\xyoption{knot}
\xyoption{curve}
\xyoption{arrow}
\xyoption{all}

\newcommand{\C}{\mathbb{C}}
\newcommand{\Q}{\mathbb{Q}}
\newcommand{\R}{\mathbb{R}}
\newcommand{\N}{\mathbb{N}}
\newcommand{\PP}{\mathbb{P}}
\newcommand{\Z}{\mathbb{Z}}

\newcommand{\Oscr}{{\mathcal O}}

\newcommand{\vtx}[1]{*+[o][F-]{\scriptscriptstyle #1}}

\newcommand{\wis}[1]{{\text{\em \usefont{OT1}{cmtt}{m}{n} #1}}}

\newtheorem{theorem}{Theorem}
\newtheorem{proposition}{Proposition}
\newtheorem{lemma}{Lemma}

{\theorembodyfont{\rmfamily} 
\newtheorem{example}{Example}

\newtheorem{definition}{Definition} }

\preprint{UIA preprint 2002-07}

\title{Isolated singularities, smooth orders \\
and Auslander regularity}

\author{Raf Bocklandt  \thanks{The author wishes to thank the university of Bielefeld for his stay in the period september-october 2002} \\
Departement Wiskunde en Informatica, Universiteit Antwerpen  \\
B-2020 Antwerp (Belgium) \\
E-mail :  \email{rafael.bocklandt@ua.ac.be}
}

\author{Lieven Le Bruyn \\
Departement Wiskunde en Informatica, Universiteit Antwerpen  \\
B-2020 Antwerp (Belgium) \\
E-mail :  \email{lieven.lebruyn@ua.ac.be}}

\author{Stijn Symens \thanks{Research assistant of the Fund for Scientific Research-Flanders (Belgium)} \\
Departement Wiskunde en Informatica, Universiteit Antwerpen  \\
B-2020 Antwerp (Belgium) \\
E-mail :  \email{stijn.symens@ua.ac.be}
}

\abstract{In this note we prove that a smooth order satisfying the reverse geometric
engineering conditions in stringtheory \cite[\S 2]{Berensteinreverse} (in any compactifying dimension) is Auslander regular. Moreover, we classify the \'etale local structure of smooth orders over an isolated central singularity.}

\begin{document}

There are at least two well-developed theories of noncommutative smoothness. Auslander regularity (generalizing Serre's homological characterization of commutative regular algebras)
is well suited to deal with questions from K-theory, derived categories and intersection theory.
Formally smooth algebras and smooth orders (generalizing Grothendieck's categorical characterization of commutative regular algebras) are well suited to deal with geometric questions such as the \'etale local structure (of the algebra and its center) and Brauer-Severi varieties.
Not surprisingly, 'real life situations' (such as stringtheory) often require the best of both worlds. 

\section{Definitions}

Let $A$ be an affine $\C$-algebra and denote with $\wis{rep}_n~A$  the affine scheme of $n$-dimensional representations of $A$. The basechange group $GL_n$ acts on this scheme and the geometric points of the algebraic quotient $\wis{iss}_n~A = \wis{rep}_n~A// GL_n$ classify the isomorphism classes of semisimple $n$-dimensional representations of $A$. In general, $\wis{rep}_n~A$ can have several connected components and in the decomposition
\[
\wis{rep}_n~A = \bigsqcup_{\alpha}~\wis{rep}_{\alpha}~A \]
we say that $\alpha$ is a dimension vector of total dimension $| \alpha | = n$. The corresponding
algebraic quotient will be denoted by $\wis{iss}_{\alpha}~A$ and its coordinate ring
$Z_{\alpha} = \C[\wis{iss}_{\alpha}~A]$ is a central subring of the algebra of $GL_n$-equivariant maps
\[
\int_{\alpha}~A = M_n(\C[\wis{rep}_{\alpha}~A])^{GL_n} \]
from $\wis{rep}_{\alpha}~A$ to $M_n(\C)$. The algebra $\int_{\alpha}~A$ is a Noetherian algebra and is a finite module over $Z_{\alpha}$. We define two $\alpha$-relative notions of noncommutative smoothness on $A$.

\begin{definition} $A$ is said to be $\alpha$-Auslander regular (or equivalently, $\int_{\alpha}~A$ is Auslander regular) if the following conditions are satisfied for $B = \int_{\alpha}~A$ :
\begin{enumerate}
\item{$B$ has finite global dimension, $gldim(B) < \infty$.}
\item{For every finitely generated left $B$-module $M$, every integer $j \geq 0$ and every (right)
$B$-submodule $N$ of $Ext^j_B(M,B)$ we have that $j(N) \geq j$, Here, $j(N)$ is the {\em grade number}
of $N$ which is the least integer $i$ such that $Ext^i_B(N,B) \not= 0$.}
\item{For every finitely generated left $B$-module $M$ we have the equality
\[
GKdim(M) + j(M) = GKdim(B) \]
where $GKdim$ denotes the Gelfand-Kirillov dimension, see for example \cite{KrauseLenegan}.}
\end{enumerate}
\end{definition}

A major application of this notion is that it allows us to study finitely generated $B$-modules in terms of pure modules using the spectral sequence
\[
E_2^{p,-q}(M) = Ext^p_B(Ext^q_B(M,B),B) \Rightarrow H^{p-q}(M) \]
where $H^0(M) = M$ and $H^i(M) = 0$ for $i \not= 0$. By property $(2)$ the second term of this sequence is triangular.

\begin{definition} $A$ is said to be $\alpha$-smooth (or equivalently, $\int_{\alpha}~A$ is a smooth order) if the following conditions are satisfied :
\begin{enumerate}
\item{ The connected component $\wis{rep}_{\alpha}~A$ is a smooth variety.}
\item{A Zariski open subset $\wis{azu}_{\alpha}~A$ (the Azumaya locus of $A$) of $\wis{rep}_{\alpha}~A$ consists of simple representations.}
\end{enumerate}
\end{definition}

By $(2)$ the quotient $\wis{rep}_{\alpha}~A \rOnto^{\pi} \wis{iss}_{\alpha}~A$ is generically a principal $PGL_n$-fibration and hence determines a central simple algebra of dimension $n^2$
(where $n = | \alpha |$) over the function field of $Z_{\alpha}$. By $(1)$, $Z_{\alpha}$ is integrally closed and therefore $\int_{\alpha}~A$ is an order in the central simple algebra having as its center $Z_{\alpha}$.
The main application of this notion is that it allows us to describe the \'etale local structure of
$\int_{\alpha}~A$ and of $Z_{\alpha}$. Let $\xi$ be a point of $\wis{iss}_{\alpha}~A$ with corresponding semi-simple $n$-dimensional representation
\[
M = S_1^{\oplus e_1} \oplus \hdots \oplus S_k^{\oplus e_k} \]
Consider the {\em local quiver setting} $(Q,\epsilon)$ where $Q$ is the finite quiver on $k$ vertices $\{ v_1,\hdots,v_k \}$ (corresponding
to the distinct irreducible components of $M$) such that the number of oriented arrows from $v_i$ to
$v_j$ is given by the dimension of the extension space $Ext^1_A(S_i,S_j)$. The dimension vector
$\epsilon$ of the quiver $Q$ is given by the multiplicities $(e_1,\hdots,e_k)$ with which these
simple components occur in $M$. To be precise, there is a $GL_n$-equivariant \'etale local isomorphism
between $\wis{rep}_{\alpha} A$ and the associated fiber bundle
\[
GL_n \times^{GL(\epsilon)} \wis{rep}_{\epsilon} Q \]
where $\wis{rep}_{\epsilon} Q$ is the vectorspace of $\epsilon$-dimensional representations of the
quiver $Q$ on which the group $GL(\epsilon) = GL_{e_1} \times \hdots \times GL_{e_k}$ acts by
basechange. Moreover, the embedding $GL(\epsilon) \rInto GL_n$ is determined by the dimensions $d_i$
of the simple components $S_i$. As a consequence, there is an \'etale local isomorphism between
$\wis{iss}_{\alpha} A$ and the quotient variety $\wis{iss}_{\epsilon} Q = \wis{rep}_{\epsilon} Q // GL(\epsilon)$,
the variety parametrizing isoclasses of semisimple $\epsilon$-dimensional representations of $Q$.
In particular this allows us to control the central singularities which were classified in low dimensions in \cite{BockLBVdW}. 

\section{Auslander regularity}

In reverse geometric engineering of singularities in
stringtheory (see e.g. \cite[\S 2]{Berensteinreverse}) one is interested in the case that the non-Azumaya locus (the 'non bulk representations' in physical lingo) consists  of isolated singularities. 
We are able to determine the \'etale local structure of such $\alpha$-smooth orders.

\begin{lemma} \label{etaleclass} Let $A$ be an $\alpha$-smooth order such that $\{ p \}$ is an isolated singularity
which is locally the non-Azumaya locus. Then, the \'etale local structure of $\int_{\alpha}~A$ in
$p$ is determined by a quiver setting
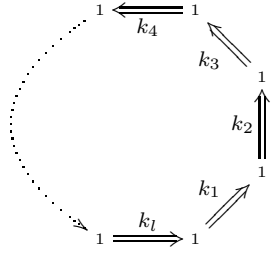
\[
\xy 0;/r.15pc/:
\POS (0,0)  *+{\txt{\tiny $1$}}="a",
(20,0) *+{\txt{\tiny $1$}}="b",
(34,14) *+{\txt{\tiny $1$}} ="c",
(34,34) *+{\txt{\tiny $1$}}="d",
(20,48) *+{\txt{\tiny $1$}} ="e",
(0,48) *+{\txt{\tiny $1$}}="f"
\POS"a" \ar@{=>}^{k_l} "b"
\POS"b" \ar@{=>}^{k_1} "c"
\POS"c" \ar@{=>}^{k_2} "d"
\POS"d" \ar@{=>}^{k_3} "e"
\POS"e" \ar@{=>}^{k_4} "f"
\POS"f" \ar@/_7ex/@{.>} "a"
\endxy 
\]
where $Q$ has $l$ vertices and all $k_i \geq 1$. The central dimension
is
\[
d = \sum_i k_i + l - 1 \]
\end{lemma}

\Proof
To start, $\epsilon$ is the dimension vector of a simple representation of $Q$. By the results of
\cite{LBProcesi} this implies that $Q$ is a strongly connected quiver (any pair of vertices $v_i$, $v_j$
is connected by an oriented path in $Q$ starting at $v_i$ and ending in $v_j$) and that the dimension vector
$\epsilon$ satisfies the numerical conditions
\[
\chi_Q(\epsilon,\delta_i) \leq 0 \qquad \text{and} \qquad \chi_Q(\delta_i,\epsilon) \leq 0 \]
where $\chi_Q$ is the Euler-form of the quiver $Q$ (that is, the bilinear form on $\Z^k$ determined by the
$k \times k$ matrix whose $(i,j)$-entry is $\delta_{ij}-$ the number of arrows from $v_i$ to $v_j$ and
where $\delta_i$ is the basevector concentrated in $v_i$). Next, we claim that $\epsilon = (1,\hdots,1)$.
If not, there are $\epsilon$-dimensional semi-simple representations of $Q$ of representation type
\[
(1,(1,\hdots,1);e_1-1,\delta_1;\hdots;e_k-1,\delta_k) \]
(the first factor indeed corresponds to a simple representation of $Q$ as $Q$ is strongly connected)
which is impossible by $GL_n$-equivariance and the fact that the non-Azumaya locus is concentrated in $p$
(which corresponds to the point of representation type $(e_1,\delta_1;\hdots;e_k,\delta_k)$).

We claim that every oriented cycle in $Q$ has as its support all the vertices $\{ v_1,\hdots,v_k \}$
and consequently that the quiver setting $(Q,\epsilon)$ is of the following shape :
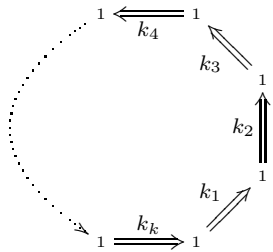
\[
\xy 0;/r.15pc/:
\POS (0,0) *+{\txt{\tiny $1$}} ="a",
(20,0) *+{\txt{\tiny $1$}} ="b",
(34,14) *+{\txt{\tiny $1$}} ="c",
(34,34) *+{\txt{\tiny $1$}} ="d",
(20,48) *+{\txt{\tiny $1$}} ="e",
(0,48) *+{\txt{\tiny $1$}} ="f"
\POS"a" \ar@{=>}^{k_k} "b"
\POS"b" \ar@{=>}^{k_1} "c"
\POS"c" \ar@{=>}^{k_2} "d"
\POS"d" \ar@{=>}^{k_3} "e"
\POS"e" \ar@{=>}^{k_4} "f"
\POS"f" \ar@/_7ex/@{.>} "a"
\endxy 
\]
Indeed, let $C$ be an oriented cycle of minimal support in $Q$, let $\delta_{iC}=1$ iff $v_i \in supp(C_)$
and zero otherwise and let $\delta_C = (\delta_{1C},\hdots,\delta_{kC})$. Then, if $C \not= \{ v_1,\hdots,v_k \}$
there would be points of representation type
\[
(1,\delta_C;e_1-\delta_{1C},\delta_1;\hdots;e_k-\delta_{kC},\delta_k) \]
contradicting the assumptions. That is, the Euler-form of the quiver $Q$ is given by the matrix
\[
\begin{bmatrix}
1 & -k_1 & 0 & \hdots & \hdots & 0 \\
0 & 1 & -k_2 & &  & 0 \\
\vdots & & \ddots & \ddots & & \vdots \\
\vdots & & & \ddots & & -k_{k-1} \\
-k_k & 0 & 0 & \hdots & \hdots & 1 
\end{bmatrix}
\]
and the statement on the dimension follows again from \cite{LBProcesi}.

\par \vskip 3mm

\begin{theorem} If $\int_{\alpha}~A$ is a smooth order such that its non-Azumaya locus consists of isolated singularities, then $\int_{\alpha}~A$ is Auslander regular.
\end{theorem}

\Proof
From \cite{LBProcesi} we recall that the ring of polynomial $GL(\epsilon)$-invariants on
$\wis{rep}_{\epsilon} Q$ is generated by the traces along oriented cycles in $Q$. Therefore,
\[
\C[\wis{iss}_{\epsilon} Q] = \C[\wis{rep}_{\epsilon} Q]^{GL(\epsilon)} = R = \C[ x_{i_1}(1)x_{i_2}(2) \hdots x_{i_k}(k) ; 1 \leq i_j \leq k_j] \subset \C[\wis{rep}_{\epsilon} Q]
\]
Therefore, $p$ is a singular point of $Z_{\alpha} = \wis{iss}_{\alpha} A$ if and only if at least two of the
$k_i \geq 2$ because by the \'etale local isomorphism the completion of the coordinate ring
$\C[ \wis{iss}_{\alpha} A]$ at the maximal ideal determined by $p$ is isomorphic to $\hat{R}$, the
completion of $R$ at the maximal ideal generated by all traces along oriented cycles in $Q$.

Similarly, we can determine $\hat{\int_{\alpha}} A$, the completion of the smooth order $\int_{\alpha} A$ at the central maximal ideal determined by $p$ from \cite{LBProcesi},
\[
\widehat{\int_{\alpha}} A \simeq
\begin{bmatrix}
M_{d_1}(\hat{R}) & M_{d_1 \times d_2}(M_{12})& M_{d_1 \times d_3}(M_{13}) & \hdots & M_{d_1 \times d_k}(M_{1k}) \\
M_{d_2 \times d_1}(M_{21}) & M_{d_2}(\hat{R}) & M_{d_2 \times d_3}(M_{23}) & \hdots & M_{d_2 \times d_k}(M_{2k}) \\
M_{d_3 \times d_1}(M_{31}) & M_{d_3 \times d_2}(M_{32}) & M_{d_3}(\hat{R}) & \hdots & M_{d_3 \times d_k}(M_{3k}) \\
\vdots & \vdots & \vdots & \ddots & \vdots \\
M_{d_k \times d_1}(M_{k1}) & M_{d_k \times d_2}(M_{k2}) & M_{d_k \times d_3}(M_{k3}) & \hdots & M_{d_k}(\hat{R})
\end{bmatrix}
\]
Here, $d_i = dim_{\C}(S_i)$ and $M_{ij}$ is the $\hat{R}$-submodule of $\C[[x_i(j),1 \leq i \leq k_j,1\leq j \leq k]]$ generated by the oriented paths in $Q$ from $v_i$ to $v_j$.

The strategy to prove that $\hat{\int_{\alpha}} A$ is Auslander regular is to use the
$GL(\epsilon) = \C^* \times \hdots \times \C^*$-action on $\wis{rep}_{\epsilon} Q$ to
obtain central elements of $\hat{\int_{\alpha}} A$ corresponding to certain arrows. Modding
out these elements in a specific order will reduce the quiver until we are left with a
hereditary (in particular, Auslander-regular) order. We can then retrace our steps using the
result that if $B$ is a Noetherian algebra with central element $c$ such that $B/(c)$ is Auslander-regular,
then so is $B$. 

Let $\{ z_1,\hdots,z_l \}$ be the vertex-indices such that $k_{z_i} \geq 2$ and after
cyclically renumbering the vertices (if needed) we may assume that $z_l = k$. For all $i \not= z_1$
set $x_i(1) = 1$ and let $\hat{R}_1$ respectively $\hat{B}_1$ be the algebra obtained from
$\hat{R}$ respectively $\hat{\int_{\alpha}} A$ using these assignments, then clearly,
\[
\hat{R} \simeq \hat{R}_1 \qquad \text{and} \qquad \widehat{\int_{\alpha}} A \simeq \hat{B}_1 \]
The advantage of this change of generators is that $\{ x_{z_1}(1),\hdots,x_{z_1}(k_{z_1}) \}$ are
generators of $\hat{R}_1$ and the quotient-algebras
\[
\overline{R}_1 = \frac{\hat{R}_1}{(x_{z_1}(2),\hdots,x_{z_1}(k_{z_1}))} \qquad \text{respectively}
\qquad \overline{B}_1 = \frac{\hat{B}_1}{(x_{z_1}(2),\hdots,x_{z_1}(k_{z_1}))}
\]
are isomorphic to the completion of the algebra of polynomial invariants (resp. equivariant maps)
of the quiver setting $(Q_1,\epsilon)$ where $Q_1$ has the same shape as $Q$ except that there
is just one arrow from $v_{z_1}$ to $v_{z_1+1}$. Repeat this procedure, starting with the quiver setting $(Q_1,\epsilon)$ with vertex $v_{z_2}$. That is, for all $i \not= z_2$ set $x_i(1)=1$ and
let $\hat{R}_2$ respectively $\hat{B}_2$ be the algebra obtained from the completion of the algebra of polynomial invariants (resp. equivariant maps)
of the quiver setting $(Q_1,\epsilon)$ using these assignments and let $\overline{R}_2$ respectively $\overline{B}_2$ be the quotient algebras obtained by modding out the
generators $\{ x_{z_2}(2),\hdots,x_{z_2}(k_{z_2}) \}$ of $\hat{R}_2$ and observe that these
quotients are the relevant algebras corresponding to a quiver setting $(Q_2,\epsilon)$ where
$Q_2$ has the same shape as $Q_1$ except that there is just one arrow from $v_{z_2}$ to
$v_{z_2+1}$ and so on.

After $l$ iterations of this procedure we arrive at the quiver setting $(Q_l,\epsilon)$ where
$Q_l$ is of the form
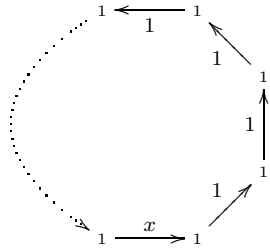
\[
\xy 0;/r.15pc/:
\POS (0,0) *+{\txt{\tiny $1$}} ="a",
(20,0) *+{\txt{\tiny $1$}} ="b",
(34,14) *+{\txt{\tiny $1$}} ="c",
(34,34) *+{\txt{\tiny $1$}} ="d",
(20,48) *+{\txt{\tiny $1$}} ="e",
(0,48) *+{\txt{\tiny $1$}} ="f"
\POS"a" \ar@{->}^{x} "b"
\POS"b" \ar@{->}^{1} "c"
\POS"c" \ar@{->}^{1} "d"
\POS"d" \ar@{->}^{1} "e"
\POS"e" \ar@{->}^{1} "f"
\POS"f" \ar@/_7ex/@{.>} "a"
\endxy 
\]
from which we deduce that $\overline{R}_l \simeq \C[[x]]$ and that
\[
\overline{B}_l \simeq \begin{bmatrix}
M_{d_1}(\C[[x]]) & M_{d_1 \times d_2}(\C[[x]])& M_{d_1 \times d_3}(\C[[x]]) & \hdots & M_{d_1 \times d_k}(\C[[x]]) \\
M_{d_2 \times d_1}(x\C[[x]]) & M_{d_2}(\C[[x]]) & M_{d_2 \times d_3}(\C[[x]]) & \hdots & M_{d_2 \times d_k}(\C[[x]]) \\
M_{d_3 \times d_1}(x\C[[x]]) & M_{d_3 \times d_2}(x\C[[x]]) & M_{d_3}(\C[[x]]) & \hdots & M_{d_3 \times d_k}(\C[[x]]) \\
\vdots & \vdots &  & \ddots & \vdots \\
M_{d_k \times d_1}(x\C[[x]]) & M_{d_k \times d_2}(x\C[[x]]) & M_{d_k \times d_3}(x\C[[x]]) & \hdots & M_{d_k}(\C[[x]])
\end{bmatrix}
\]
It is well known that $\overline{B}_l$ is an Auslander-regular algebra and as we divided out
central elements in each step (and in each step, the localizations at these central elements are Azumaya algebras with regular center hence Auslander regular), we derive using \cite[theorem III.3.6]{LiHuishiFVO} that also $\hat{\int_{\alpha}} A$
is Auslander-regular. Because Auslander-regularity is preserved under central \'etale extensions
and because $A$ is at all other points an Azumaya algebra over a commutative regular ring,
Auslander regularity of $A$ follows.

\section{Isolated singularities}

In this section we will give the \'etale local structure of an $\alpha$-smooth order in an isolated central singularity. That is, we will extend lemma~\ref{etaleclass} without the condition on the Azumaya locus. Because an $\alpha$-smooth order is locally determined by a quiver setting (see section 1) the problem reduces to classifying all quiver settings $(Q,\alpha)$ such that $\wis{iss}_{\alpha}~Q$ is an isolated singularity. 

If $v$ is a vertex having no loops in $Q$ and such that $\chi_Q(\epsilon_v,\alpha) \geq 0$ or $\chi_Q(\alpha,\epsilon_v) \geq 0$, then we replace the quiver setting $(Q,\alpha)$ by $(Q',\alpha')$ where $Q'$ is the quiver obtained from $Q$ by deleting the vertex $v$ and adding arrows corresponding to $2$-paths through $v$
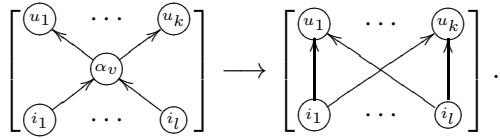
\[
\left[ \vcenter{
\xymatrix@=.3cm{
\vtx{u_1}&\cdots &\vtx{u_k}\\
&\vtx{\alpha_v}\ar[ul]\ar[ur]&\\
\vtx{i_1}\ar[ur]&\cdots &\vtx{i_l}\ar[ul]}}
\right]
\longrightarrow
\left[\vcenter{
\xymatrix@=.3cm{
\vtx{u_1}&\cdots &\vtx{u_k}\\
&&\\
\vtx{i_1}\ar[uu]\ar[uurr]&\cdots &\vtx{i_l}\ar[uu]\ar[uull]}}
\right].
\]
(note that some of the vertices in the picture may coincide leading to loops). The dimension vector $\alpha' = \alpha \mid \wis{supp}~Q'$. The reduction step $(Q,\alpha) \rTo (Q',\alpha')$ will be denoted by $R_I^v$.

\begin{theorem} Let $A$ be an $\alpha$-smooth order and $p$ a central isolated singularity. Then, the \'etale local structure of $\int_{\alpha}~A$ in $p$ is determined by a quiver setting $(Q,\epsilon)$ which can be reduced, via iterated use of $R_I^v$, to a quiver setting
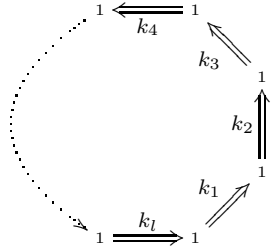
\[
\xy 0;/r.15pc/:
\POS (0,0)  *+{\txt{\tiny $1$}}="a",
(20,0) *+{\txt{\tiny $1$}}="b",
(34,14) *+{\txt{\tiny $1$}} ="c",
(34,34) *+{\txt{\tiny $1$}}="d",
(20,48) *+{\txt{\tiny $1$}} ="e",
(0,48) *+{\txt{\tiny $1$}}="f"
\POS"a" \ar@{=>}^{k_l} "b"
\POS"b" \ar@{=>}^{k_1} "c"
\POS"c" \ar@{=>}^{k_2} "d"
\POS"d" \ar@{=>}^{k_3} "e"
\POS"e" \ar@{=>}^{k_4} "f"
\POS"f" \ar@/_7ex/@{.>} "a"
\endxy 
\]
with $l \geq 2$ vertices and all $k_i \geq 2$. The central dimension
is equal to
\[
d = \sum_i k_i + l - 1 \]
\end{theorem}

Contrary to the situation of the previous section, we can have central points $q \in \wis{iss}_{\alpha}~A$ corresponding to proper semi-simple representations
\[
M = S_1^{\oplus e_1} \oplus \hdots \oplus S_l^{\oplus e_l} \]
such that $\wis{iss}_{\alpha}~A$ is smooth in $q$, or equivalently, that the {\em local quiver setting}
$(Q_q,\epsilon_q)$ defined in section 1 is {\em coregular}, that is, $\wis{iss}_{\epsilon_q}~Q_q$ is a smooth variety. Thanks to \cite{Bocklandt2002} we have a classification of coregular quiver settings. For a quiver setting $(Q,\alpha)$ with a vertex $v$ such that $\alpha_v = 1$ and there are loops in $v$ we define the reduction step $R_{II}^v$ to be $(Q,\alpha) \rTo (Q',\alpha)$ where $Q'$ is the quiver obtained from $Q$ by removing the loops in $v$.

For a quiver setting $(Q,\alpha)$ and a vertex $v$ such that $\alpha_v = k > 1$, there is a unique loop in $v$ and the neighborhood of $Q$ in $v$ is one of the situations on the left hand side of the pictures below

\vspace{1cm}

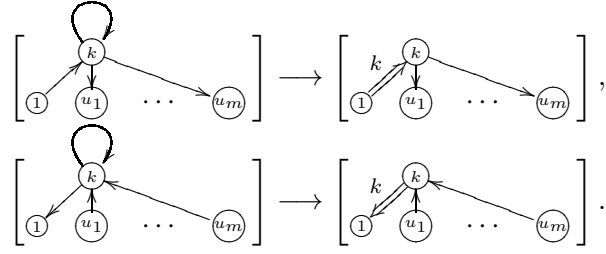
\[
\left[\vcenter{
\xymatrix@=.3cm{
&\vtx{k}\ar[d]\ar[drr]\ar@(lu,ru)&&\\
\vtx{1}\ar[ur]&\vtx{u_1}&\cdots &\vtx{u_m}}}
\right]
\longrightarrow
\left[\vcenter{
\xymatrix@=.3cm{
&\vtx{k}\ar[d]\ar[drr]&&\\
\vtx{1}\ar@2[ur]^k&\vtx{u_1}&\cdots &\vtx{u_m}}}
\right],
\]

\vspace{.2cm}

\[
\left[\vcenter{
\xymatrix@=.3cm{
&\vtx{k}\ar[dl]\ar@(lu,ru)&&\\
\vtx{1}&\vtx{u_1}\ar[u]&\cdots &\vtx{u_m}\ar[ull]}}
\right]
\longrightarrow
\left[\vcenter{
\xymatrix@=.3cm{
&\vtx{k}\ar@2[dl]_k&&\\
\vtx{1}&\vtx{u_1}\ar[u]&\cdots &\vtx{u_m}\ar[ull]}}
\right].
\]
(again, some of the vertices may be the same). Then we define a reduction step
$R_{III}^v$ which sends  $(Q,\alpha) \rTo (Q',\alpha)$ where $Q'$ is $Q$ with the neighborhood of $v$ replaced by the situation on the right hand side of the pictures. The main result of \cite{Bocklandt2002} asserts that $(Q,\alpha)$ is a coregular quiver setting if and only if it can be reduced by an iterated use of the reduction steps $R_I^v,R_{II}^v$ and $R^v_{III}$ (and their inverses) to one of the three quiver settings below :
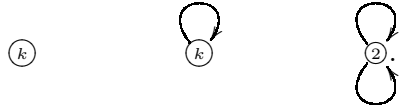
\vspace{.5cm}
\[
\xymatrix{\vtx{k}}\hspace{2cm} \xymatrix{\vtx{k}\ar@(lu,ru)}\hspace{2cm} \xymatrix{\vtx{2}\ar@(lu,ru)\ar@(ld,rd)}.
\]

\vspace{.5cm}

A {\em representation type} $\tau = (e_1,\beta_1;\hdots;e_l,\beta_l)$ of a quiver setting $(Q,\alpha)$ satisfies $\alpha = e_1 \beta_1 + \hdots + e_l \beta_l$ and all $\beta_i$ are dimension vectors of simple representations of $Q$ (and we have a description of those from
\cite{LBProcesi}). The local quiver setting in a point $\xi \in \wis{iss}_{\alpha}~Q$ of representation type $\tau$ depends only on $\tau$ : $Q_{\tau}$ is the quiver on $l$ vertices such that there are exactly $\delta_{ij} - \chi_Q(\beta_i,\beta_j)$ arrows (or loops) from the $i$-th to the $j$-th vertex and $\alpha_{\tau} = (e_1,\hdots,e_l)$, see \cite{LBProcesi}. The {\em stratum} $S_{\tau}$ consisting of all points in $\wis{iss}_{\alpha}~Q$ having representation type $\tau$ has dimension
\[
dim~S_{\tau} = \sum_{loop} (l_j-1)e_j^2 + 1 \]
where the sum is taken over all vertices $w_j$ having loops in $Q_{\tau}$, see \cite{LBProcesi}. If we apply this to the representation type $(\alpha_1,\epsilon_1;\hdots;\alpha_k,\epsilon_k)$ of the trivial representation we deduce :

\begin{lemma} If $(Q,\alpha)$ is a quiver setting such that $\wis{iss}_{\alpha}~Q$ is an isolated singularity, then there are no loops in $Q$.
\end{lemma}

If $(Q,\alpha) \rDotsto (Q',\alpha')$ is a sequence of reductions $R_I^v,R_{II}^v$ or $R_{III}^v$ we have that either
\[
\wis{iss}_{\alpha}~Q = \wis{iss}_{\alpha'}~Q' \qquad \text{or} \qquad \wis{iss}_{\alpha}~Q = \wis{iss}_{\alpha'}~Q' \times \C^z \]
for some $z$ (see \cite{Bocklandt2002}). By this and the lemma we have that any reduction of a quiver setting $(Q,\alpha)$ with $\wis{iss}_{\alpha}~Q$ an isolated singularity involves only reduction steps $R_I^v$. We will characterize the {\em reduced} settings, that is those that cannot be reduced further.

\begin{lemma} If $(Q,\alpha)$ is a reduced quiver setting with $\wis{iss}_{\alpha}~Q$ an isolated singularity, then $\alpha = \mathbf{1} = (1,\hdots,1)$.
\end{lemma}

\Proof
Assume $v$ is a vertex having maximal $\alpha_v \geq 2$. Because $(Q,\alpha)$ is reduced it follows from the definition of reduction step $R_I^w$ that for all vertices $w$ we have
\[
\chi_Q(\epsilon_w,\alpha) < 0 \quad \text{and} \quad \chi_Q(\alpha,\epsilon_w) < 0 \]
Therefore, by $\cite{LBProcesi}$ we have that $\alpha - \epsilon_v$ is the dimension vector of a simple representation of $Q$ and we look at the local quiver setting $(Q_{\tau},\alpha_{\tau})$ for the representation type $\tau = (1,\epsilon_v;1,\alpha-\epsilon_v)$. This is of the form
$$
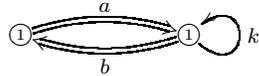
\xymatrix{\vtx{1} \ar@/^/@2{->}[rr]^a &  & \vtx{1} \ar@(rd, ru)[]_k \ar@/^/@2{->}[ll]^b }
$$
where $a = -\chi_Q(\epsilon_v,\alpha-\epsilon_v) = -\chi_Q(\epsilon_v,\alpha)+1 \geq 2$, $b = -\chi_Q(\alpha-\epsilon_v,\epsilon_v) = -\chi_Q(\alpha,\epsilon_v) + 1 \geq 2$ and
$k = 1 -\chi_Q(\alpha-\epsilon_v,\alpha-\epsilon_v)$. Therefore, $(Q_{\tau},\alpha_{\tau})$ is not coregular and hence $\wis{iss}_{\alpha}~Q$ is not an isolated singularity, a contradiction.

\par \vskip 6mm

\noindent
{\bf Proof of theorem 2 : } If $(Q,\alpha)$ is a quiver setting such that $\wis{iss}_{\alpha}~Q$ is an isolated singularity, we can reduce it by iterated use of $R_I^v$ to a setting $(Q',\mathbf{1})$ by the previous lemma. We now claim that $Q'$ is of the prescribed form, that is that every cycle in $Q'$ passes through all vertices. If not let
$\{ v_{i1},\hdots,v_{ip} \}$ be the vertices through which a cycle does {\em not}  pass and consider the representation type
\[
\tau = (1,\mathbf{1}-\epsilon_{v_{i1}}- \hdots-\epsilon_{v_{ip}};1,\epsilon_{v_{i1}};\hdots;1,\epsilon_{v_{ip}}) \]
The local quiver $Q_{\tau}$ has $p+1$ vertices $\{ w_0,w_1,\hdots,w_p \}$ where $w_j$ corresponds to $v_{ij}$ and $w_0$ collects the remaining vertices $V$ of $Q$. Via this identification, the quiver on $\{ w_1,\hdots,w_p \}$ is identical to that of $Q$ on $\{ v_{i1},\hdots,v_{ip} \}$ and the number of arrows from (resp. to) $w_0$ to (resp. from) $w_j$ is equal to the number of arrows from (resp. to) $V$ to (resp. from) $v_{ij}$ in $Q$ and there is a number of loops in $w_0$.

If $(Q_{\tau},\mathbf{1}_{\tau})$ is reduced (after removing the loops at $w_0$), then it cannot be coregular by \cite{Bocklandt2002} as $Q_{\tau}$ has at least two vertices whence $\wis{iss}_{\alpha}~Q$ is {\em not} an isolated singularity.

If $(Q_{\tau},\mathbf{1}_{\tau})$ can be reduced (after removing the loops at $w_0$), then the only possible reduction step is $R_I^{w_0}$ as all $w_j$ ($j \geq 1$) have at least two incoming and two outgoing arrows. The following lemma asserts that $(Q_{\tau},\mathbf{1}_{\tau})$ cannot be regular whence again $\wis{iss}_{\alpha}~Q$ cannot be an isolated singularity.

\begin{lemma} A coregular quiver setting $(Q,\mathbf{1})$ with $Q$ strongly connected and having more than one vertex has at least two vertices $v$ allowing reduction step $R_I^v$.
\end{lemma}

\Proof
By induction on the number $n$ of vertices. If $n=2$, then by the classification of coregular quiver settings $(Q,\mathbf{1})$ must have the form
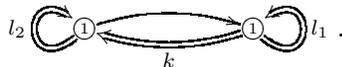
\[
\xymatrix{\vtx{1} \ar@/^/@{->}[rr]\ar@2@(ld, lu)^{l_2}  &  & \vtx{1}
\ar@2@(rd, ru)[]_{l_1} \ar@/^/@2{->}[ll]^k}.
\]
whence (after removing the loops) both vertices allow reduction $R_I$. If $n > 2$ perform one reduction $R_I^v$ (say with one outgoing arrow ending in $w$) to produce of quiver $Q'$ on $n-1$ vertices.
In $Q'$ only $w$ can change its (ir)reducible status (either way). As we removed the reducible vertex $v$ from $Q$ the number of reducible vertices in $Q'$ is less than or equal the number of reducible vertices of $Q$ but by induction $Q'$ has at least two reducible vertices.

\begin{theorem} The isolated central singularities of an $\alpha$-smooth order $A$ and a $\beta$-smooth order $B$ are \'etale equivalent if and only if their local quivers can be reduced to quiver settings 
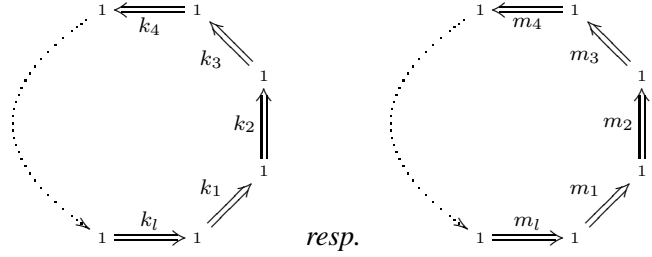
\[
\xy 0;/r.15pc/:
\POS (0,0)  *+{\txt{\tiny $1$}}="a",
(20,0) *+{\txt{\tiny $1$}}="b",
(34,14) *+{\txt{\tiny $1$}} ="c",
(34,34) *+{\txt{\tiny $1$}}="d",
(20,48) *+{\txt{\tiny $1$}} ="e",
(0,48) *+{\txt{\tiny $1$}}="f"
\POS"a" \ar@{=>}^{k_l} "b"
\POS"b" \ar@{=>}^{k_1} "c"
\POS"c" \ar@{=>}^{k_2} "d"
\POS"d" \ar@{=>}^{k_3} "e"
\POS"e" \ar@{=>}^{k_4} "f"
\POS"f" \ar@/_7ex/@{.>} "a"
\endxy 
\quad \text{resp.} \quad
\xy 0;/r.15pc/:
\POS (0,0)  *+{\txt{\tiny $1$}}="a",
(20,0) *+{\txt{\tiny $1$}}="b",
(34,14) *+{\txt{\tiny $1$}} ="c",
(34,34) *+{\txt{\tiny $1$}}="d",
(20,48) *+{\txt{\tiny $1$}} ="e",
(0,48) *+{\txt{\tiny $1$}}="f"
\POS"a" \ar@{=>}^{m_l} "b"
\POS"b" \ar@{=>}^{m_1} "c"
\POS"c" \ar@{=>}^{m_2} "d"
\POS"d" \ar@{=>}^{m_3} "e"
\POS"e" \ar@{=>}^{m_4} "f"
\POS"f" \ar@/_7ex/@{.>} "a"
\endxy 
\]
having the same number of vertices and such that the $l$-tuples $(k_1,\hdots,k_l)$ and $(m_1,\hdots,m_l)$ are the same upto a permutation.
\end{theorem}

\Proof
We give two proofs of this result. By \cite{LBProcesi} the coordinate ring of a quotient variety $\wis{iss}_{\epsilon}~Q$ is generated by traces along oriented cycles in the quiver $Q$. In the case of the left hand quiver settings, these invariants are easy to determine : the dimension of the power $\mathfrak{m}^i$ of the maximal graded ideal $\mathfrak{m}^i/\mathfrak{m}^{i+1}$ is equal to
\[
M_i = \binom{k_1+i-1}{i} \hdots \binom{k_n+i-1}{i} \]
whence $M_{i+1}/M_i = (i+1)^{-n}(i+k_1) \hdots (i+k_n)$. The rational function
\[
f(x) = \frac{(x+k_1)\hdots (x+k_n)}{(x+1)^n} \]
is determined by its values on all $x \in \N$ whence the dimension-sequence $M_{i+1}/M_i$ determines the roots and their multiplicity (note that none of the $k_i = 1$ as the quiver setting is reduced. From this the difficult part of the result follows.
Alternatively, the result follows from the fact that $\wis{iss}_{\epsilon}~Q$ is the cone on the projective variety $\PP^{k_1-1} \times \hdots \times \PP^{k_n-1}$.

\section{Applications}

A quiver gauge theory consists of a quiver setting $(Q,\alpha)$ together with the choice of a
necklace as in \cite{LBBocklandt} (a
{\em superpotential}) $W \in \wis{dR}^0_V~\C Q = \C Q / [\C Q,\C Q]$, that is,
\[
W = \sum_j a_{i_1} \hdots a_{i_{l_j}} \]
is a sum of oriented cycles in the quiver $Q$ with arrows say $\{ a_1,\hdots,a_l \}$. Such a necklace induces a $GL(\alpha)$ invariant polynomial function
\[
W~:~\wis{rep}_{\alpha}~Q \rTo \C \qquad V \mapsto \sum_j Tr(V_{a_{i_1}} \hdots V_{a_{i_{l_j}}})
\]
and hence a $GL_n$-invariant function on the $\alpha$-component $GL_n \times^{GL(\alpha)} \wis{rep}_{\alpha}~Q$ of $\wis{rep}_n~\C Q$. From the
differential
$\wis{dR}^0_V~\C Q \rTo^d \wis{dR}^1_V~\C Q$ we can define by \cite{LBBocklandt}  {\em partial
differential
operators} associated to any arrow $a$ in $Q$ with start vertex $v_i$ and end vertex $v_j$
\[
\frac{\partial}{\partial a}~:~\wis{dR}^0_{V}~\C Q \rTo e_i \C Q e_j \qquad
\text{by} \qquad df = \sum_{a \in Q_a}~\frac{\partial f}{\partial a} d a \]
To take the partial derivative of a necklace word $w$ with respect to an
arrow $a$, we run through $w$ and each time we encounter $a$ we open the
necklace by removing
that occurrence of $a$ and then take the sum of all the paths obtained.

\begin{definition} The {\em vacualgebra} of a quiver gauge theory determined by the quiver setting
$(Q,\alpha)$ and the superpotential $W$ is the Noetherian affine $\C$-algebra
\[
\int_{\alpha}~\partial_Q~W \qquad \text{where} \qquad \partial_Q~W =  \frac{\C Q}{(\frac{\partial W}{\partial a_1},\hdots,
\frac{\partial W}{\partial a_l})} \]
The affine variety $\wis{rep}_{\alpha}~\partial_Q~W$ is said to be the space of {\em vacua} and the algebraic quotient $\wis{iss}_{\alpha}~\partial_Q~W$ is called the {\em moduli of superpotential vacua}, see for example \cite{LutyTaylor}.
\end{definition}

In order to get realistic models, one has to impose additional conditions, for example that the
superpotential $W$ is cubic (meaning that every arrow in $Q$ must belong to at least one oriented cycle of length $\leq 3$) or that  $\wis{iss}_{\alpha}~\partial_Q~W$ is three-dimensional, see for
example \cite{Berenstein1}, \cite{DouglasGreeneMorrison}.

Applications of Auslander regularity of $\int_{\alpha}~\partial_Q~W$ are well-documented in the literature (a.o. \cite[\S 3]{Berenstein1},\cite[\S 6]{Berenstein2} or \cite{Berenstein3}). Applications
of smoothness of $\int_{\alpha}~\partial_Q~W$ (that is, that the space of vacua $\wis{rep}_{\alpha}~\partial_Q~W$ is a smooth variety) are more implicit. In comparing the algebraic quotient with the moment map description (comparing $F$-terms to $D$-terms) or
defining a K\"ahler metric one can get by using the induced properties from smoothness of
$\wis{rep}_{\alpha}~\C Q$. However, in comparing geometrical properties (such as flips and flops)
of related moduli spaces
(of semi-stable representations for algebraists, adding Fayet-Iliopoulos terms for physicists) one sometimes uses the stronger results of \cite{Thaddeus} and \cite{Dolgatchev} for which smoothness of the total space is crucial. Ideally, one would like to have vacualgebras having both smoothness conditions.

\begin{example}~(The conifold algebra, see \cite{Berenstein4}) The relevant quiver setting
$(Q,\alpha)$ is of the form
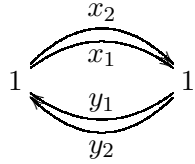
\[
\xymatrix{1 \ar@/^3ex/[rr]_{\txt{$x_1$}} \ar@/^4ex/[rr]^{\txt{$x_2$}} && 1
\ar@/^3ex/[ll]_{\txt{$y_1$}} \ar@/^4ex/[ll]^{\txt{$y_2$}}} \]
and the necklace (superpotential) is taken to be
\[
W = \lambda((x_1y_2-x_2y_1)^2-(y_1x_2-y_2x_1)^2) \]
Therefore, the defining equations of $\partial_Q~W$ are (taking into account that $x_ix_j=0$ and
$y_iy_j = 0$)
\begin{eqnarray*}
\frac{\partial W}{\partial x_1} &= y_1x_2y_2-y_2x_2y_1= 0  \\
\frac{\partial W}{\partial x_2} &= y_2x_1y_1-y_1x_1y_2 = 0 \\
\frac{\partial W}{\partial y_1} &= x_2y_2x_1-x_1y_2x_2 = 0 \\
\frac{\partial W}{\partial y_2} &= x_1y_1x_2-x_2y_1x_1= 0
\end{eqnarray*}
Observe that these identities are satisfied for {\em all} representations of $\wis{rep}_{\alpha}~Q$
as $\alpha = (1,1)$ and therefore
\[
\wis{rep}_{\alpha}~\partial_Q~W = GL_2 \times^{\C^* \times \C^*} \wis{rep}_{\alpha}~Q \quad
\text{and} \quad \int_{\alpha}~\partial_Q~W = \int_{\alpha}~\C Q \]
and therefore the vacualgebra is a smooth order.

Moreover, the quotient variety $\wis{iss}_{\alpha}~Q$ is easily seen to be the conifold singularity as the ring of invariants is generated by the primitive oriented cycles
\[
x = x_1y_1 \quad y = x_2y_2 \quad u = x_1y_2 \quad v = x_2y_1 \]
which satisfy the relation $xy = uv$. Therefore, by the theorem the {\em conifold algebra}
$\int_{\alpha}~\partial_W~Q$ is also Auslander regular. One can also check immediately that the 
description of the conifold algebra given in \cite[\S 1]{Berenstein4} is the algebra of equivariant maps from $\wis{rep}_{\alpha}~Q$ to $M_2(\C)$ (or to be more precise, if we take the relevant gauge groups into account, a ring Morita equivalent to the conifold algebra).
\end{example}

\par \vskip 3mm
By the results of \cite{BockLBVdW} we know that the only type of singularity that can occur in the center of a smooth order (in dimension three) is the conifold singularity. Therefore, in most models considered by physicists, see a.o. \cite{Berenstein1}, \cite{Greene} or \cite{Sardo} the space of superpotential vacua $\wis{rep}_{\alpha}~\partial_Q~W$ must contain singularities as the moduli space is a three dimensional quotient variety (different from the conifold) or has a one-dimensional family of singularities (which cannot happen for a three dimensional smooth order).

There is a standard way to remove (most of) the singularities in $\wis{rep}_{\alpha}~\partial_Q~W$
by restricting to semistable representations. Let us quickly run through the process. For a dimension vector $\alpha = (v_1,\hdots,v_k) \in \N^k$ let $U(\alpha)$ be the quotient of the
Lie group $U(v_1) \times \hdots \times U(v_k)$ by the one-dimensional central subgroup
$U(1)(1_{v_1},\hdots,1_{v_k})$. The {\em real moment map} for $\alpha$-dimensional
quiver representations of $Q$ is the map
\[
\wis{rep}_{\alpha}~Q \rTo^{\mu_{\R}} \wis{Lie}~U(\alpha) \qquad V \mapsto \frac{i}{2} \sum_{j=1}^l
[V_{a_j},V_{a_j}^{\dagger}] \]
There is a natural one-to-one correspondence (actually a homeomorphism)
\[
\wis{iss}_{\alpha}~Q \leftrightarrow \mu^{-1}_{\R}(0)/U(\alpha) \]
Let $\mu = (u_1,\hdots,u_k) \in \Q^k$ such that $\mu.\alpha = \sum u_iv_i = 0$, then we say that a representation $V \in \wis{rep}_{\alpha}~Q$ is $\mu$-semistable (resp.
 $\mu$-stable) if for all proper subrepresentations $W \subset V$ we have that $\mu.\beta \geq 0$
(resp. $\mu.\beta > 0$) where $\beta$ is the dimension vector of $W$. If $\wis{rep}_{\alpha}^{\mu}~Q$ denotes the Zariski open set (possibly empty) of $\mu$-semistable representations of $\wis{rep}_{\alpha}~Q$, then the geometric invariant quotient
\[
\wis{moduli}^{\mu}_{\alpha}~Q = \wis{rep}_{\alpha}^{\mu}~Q // GL(\alpha) \rOnto \wis{iss}_{\alpha}~Q  \]
classifies the isomorphism classes of $\alpha$-dimensional direct sums of $\mu$-stable representations and is a projective bundle over $\wis{iss}_{\alpha}~Q$. Moreover, there is a moment map description of this moduli space
\[
\wis{moduli}^{\mu}_{\alpha}~Q = \mu^{-1}_{\R}(\mu) / U(\alpha) \]
For more details on these matters we refer to \cite{King}.

If $I_W$ denotes the set of zeroes of the ideal of relations of $\C[\wis{rep}_{\alpha}~Q]$ imposed by the defining relations of $\partial_Q~W$, then
\[
\wis{iss}_{\alpha}~\partial_Q~W = (\mu^{-1}_{\R}(0) \cap I_W) / U(\alpha) \]
and the geometric invariant quotient of the open set of $\mu$-stable representations of $\partial_Q~W$ is a projective bundle over it
\[
\wis{moduli}^{\mu}_{\alpha}~\partial_Q~W = (\mu_{\R}^{-1}(\mu) \cap I_W) / U(\alpha) \rOnto
\wis{iss}_{\alpha}~\partial_Q~W \]
The moduli spaces of $\mu$-semistable representations can be covered by open sets determined by determinantal semi-invariants and consequently we can define a sheaf of noncommutative orders
\[
\Oscr^{\mu}_{\partial_Q~W} \qquad \text{over} \quad \wis{moduli}^{\mu}_{\alpha}~\partial_Q~W
\]
which locally is isomorphic to $\wis{iss}_{\alpha}~A$ of a suitable algebra $A$.

In favorable situations, $\wis{rep}_{\alpha}^{\mu}~\partial_Q~W$ will be a smooth variety and the moduli space $\wis{moduli}^{\mu}_{\alpha}~\partial_Q~W$ will be a (partial) desingularization
of $\wis{iss}_{\alpha}~\partial_Q~W$. In physical terminology this process is described as 'adding a Fayet-Iliopoulos term'. An immediate consequence of the theorem then implies :

\begin{proposition} With notations as before, if $\wis{rep}_{\alpha}^{\mu}~\partial_Q~W$ is a smooth variety and if the partial desingularization $\wis{moduli}_{\alpha}^{\mu}~\partial_Q~W$ has  isolated singularities as its non-Azumaya locus, then $\Oscr^{\mu}_{\partial_Q~W}$ is a sheaf of Auslander regular orders over $\wis{moduli}^{\mu}_{\alpha}~\partial_Q~W$.
\end{proposition}

In physical relevant settings, the resolution $\wis{moduli}^{\mu}_{\alpha}~\partial_Q~W \rOnto
\wis{iss}_{\alpha}~\partial_Q~W$ will often be crepant meaning that the moduli space is a Calabi-Yau manifold and the sheaf of orders will be an Azumaya sheaf. However, there may be relevant situations where we only have a partial desingularization, the remaining singularities are necessarily of conifold type and the sheaf of orders is locally Morita equivalent to the conifold algebra in the singularities. This explains the importance of conifold transitions in (partial) resolutions of three dimensional quotient singularities.

\providecommand{\href}[2]{#2}
\begingroup\raggedright\endgroup

\end{document}